\documentclass[letterpaper]{article}
\usepackage{aaai}
\usepackage{times}
\usepackage{helvet}
\usepackage{courier}
\usepackage[ruled,vlined,linesnumbered]{algorithm2e}
\usepackage{amsmath, amsfonts, amssymb, graphicx, tabularx, accents, subfigure}
\usepackage{graphicx}
\usepackage{nomencl}

\let\chapter\section

\newtheorem{thm}{Theorem}[section]

\newtheorem{obs}[thm]{Observation}

\newcommand{\SubProblem}{{\cal{Q}}(s)}

\newcommand{\Nodes}{{\cal{N}}}
\newcommand{\Edges}{{\cal{E}}}
\newcommand{\Disasters}{{\cal{S}}}
\newcommand{\Damages}{{\cal{D}}}
\newcommand{\HardenedDamages}{{\cal{D}^\prime}}

\newcommand{\lineVariable}{x}
\newcommand{\switchVariable}{\tau}
\newcommand{\hardenVariable}{t}
\newcommand{\microgridVariable}{z}
\newcommand{\facilityVariable}{u}
\newcommand{\ChanceVariable}{v}

\newcommand{\lineuseVariable}{x^s}
\newcommand{\switchuseVariable}{\tau^s}
\newcommand{\lineCycleVariable}{\bar{x}^s}
\newcommand{\switchCycleVariable}{\bar{\tau}^s}
\newcommand{\linehardenVariable}{t^s}
\newcommand{\microgriduseVariable}{z^s}
\newcommand{\facilityuseVariable}{u^s}
\newcommand{\generatorVariable}{g^s}
\newcommand{\loadVariable}{l^s}
\newcommand{\loadServeVariable}{y^s}
\newcommand{\flowVariable}{f^s}

\newcommand{\lineCost}{c}
\newcommand{\switchCost}{\kappa}
\newcommand{\hardenCost}{\psi}
\newcommand{\facilityCost}{\alpha}
\newcommand{\microgridCost}{\zeta}

\newcommand{\capacity}{Q}
\newcommand{\phases}{{\cal{P}}}
\newcommand{\phaseVariation}{\beta}
\newcommand{\demand}{d}
\newcommand{\availablecapacity}{G}
\newcommand{\microgridcapacity}{Z}
\newcommand{\cycles}{\cal{C}}
\newcommand{\criticalLoad}{\lambda}
\newcommand{\load}{\gamma}
\newcommand{\Load}{{\cal{L}}}

\immediate\write18{makeindex \jobname.nlo -s ./nomencl.ist -o \jobname.nls}
\makenomenclature
\begin{document}
%
\title{Designing Resilient Electrical Distribution Grids}
\author{Emre Yamangil$^1$, Russell Bent$^2$ and Scott Backhaus$^2$\\
$^1$ Rutgers University, New Brunswick, NJ, 08901\\
$^2$ Los Alamos National Laboratory, Los Alamos, NM, 87545
}
\maketitle
\begin{abstract}

\begin{quote}
Modern society is critically dependent on the services provided by engineered infrastructure networks.  When natural disasters (e.g. Hurricane Sandy) occur, the ability of these networks to provide service is often degraded because of physical damage to network components. One of the most critical of these networks is electric power, with medium voltage distribution circuits often suffering the  most severe damage. However, well-placed upgrades to these distribution grids can greatly improve post-event network performance.  We formulate an optimal electrical distribution grid design problem as a two-stage, stochastic mixed-integer program with damage scenarios from natural disasters modeled as a set of stochastic events. We develop and investigate the tractability of an exact and several heuristic algorithms based on decompositions that are hybrids of techniques developed by the AI and operations research communities. We provide computational evidence that these algorithms have significant benefits when compared with commercial, mixed-integer programming software.

\end{quote}
\end{abstract}

\section{Introduction}
\label{sec:intro}

Natural disasters such as earthquakes, hurricanes, and other extreme weather pose serious risks  to modern critical infrastructure including electrical distribution grids. At the peak of Hurricane Sandy, 65\% of New Jersey's customers lost power \cite{Sandy:13}. Recent U.S. government sources \cite{WhiteHouse2013,DOEResilience2013} suggest that new methodologies for improving system resilience to these events is necessary.
Here, we focus on developing methods for designing and upgrading distribution grids to better withstand and recover from these threats that are inspired by techniques developed in the artificial intelligence and operations research communities. Our approach minimizes the upgrade budget while meeting a minimum standard of service by selecting from a set of potential upgrades,  e.g. adding redundant lines, adding distributed (microgrid) generation (i.e. wind,  solar, and combined heat and power), hardening existing components, etc.

We formulate our approach, i.e. Optimal Resilient Distribution Grid Design (ORDGDP), as a two-stage mixed-integer program. The first (investment) stage selects from the set of potential upgrades to the network.  The second (operations) stage evaluates the network performance benefit of the upgrades against a set of damage scenarios sampled from a stochastic distribution.
We first develop an exact solution method that shares a number of similarities with Benders Decomposition \cite{VW:10} by exploiting decomposition across the sampled scenarios.  We also develop a metaheuristic that we call Scenario Based Variable Neighborhood Decomposition Search (SBVNDS) that is a hybrid of Variable Neighborhood Search \cite{LHMU:10} and the exact method.
We present numerical evidence that our exact method is more efficient than out-of-the-box commercial mixed-integer programming solvers, and that our heuristic achieves near-optimal results in a fraction of the time required by exact methods.

\paragraph{Literature Review}  Network design problems and their variations are generally NP-complete \cite{TPZ:10,NPTZ:13,Johnson1978}. However, recent work by \cite{BBT:10} demonstrates that AI-based methods can lead to substantial improvement for realistic applications.

While the specific problem of designing resilient distribution systems is novel, a number of related problems exist. The flow of electric power in tree-like distribution networks is related to multi-commodity network flows making our problem similar to the design of multi-commodity flow networks with stochastic link and edge failures \cite{SAGS:03,GS:08}. However, the second stage of our formulation requires binary variables making our problem considerably more difficult than typical second-stage flow problems.
The interdiction literature includes related max-min or min-max problems where the goal is to operate or design a system to make it as resilient as possible to an adversary who can damage up to $k$ elements. Such models are similar to ours if a $k$ is chosen that bounds the worst-case disaster \cite{Chen:13,Chen2014a,Salmeron2009,Delgadillo2010}.
Binary variables at all stages make these models computationally challenging and solvable only for small $k$. Here, we exploit the probabilistic nature of our adversary to increase the size of tractable problems (eliminates a stage of binary variables).

In power engineering, papers have primarily focused on resilient system operation \cite{GFW:14,Li2014,Khushalani2007} using controls such as line switching. The ORDGDP  is a fundamental generalization of the resilient operations problem because 1) this problem is embedded in our second stage and 2) minimizing the number of switch actions \cite{Li2014} can be thought of as a design problem for a single scenario.  Finally, there is also a general power grid expansion planning problem for stochastic events \cite{Jabr2013} that is a variant of the the single commodity flow problem, with the twist that flows are not directly controllable. Like stochastic multi-commodity flow, the second-stage variables are not binary.

The key contributions of the paper include:

\begin{itemize}
\item Computationally efficient algorithms for solving stochastic network design problems with discrete variables at each stage. The algorithms are based on hybrid optimization methods similar to recent work that combines Bender's Decomposition with heuristic master solutions  \cite{RBH:14}.

\item Introduction of a problem of critical importance to energy problems where AI researchers can make significant contributions. AI has made many recent significant contributions to energy problems \cite{HGC:12,RV:12,CHB:12,GJN:13,RV:13,JNN:14,RV:11,SS:13,TCHS:13}.

\end{itemize}


\section{Problem Description}
\label{sec:desc}

%

\nomenclature[AA]{$\textbf{Parameters}$}{}%
\nomenclature[AB]{$\Nodes$}{set of nodes (buses).}%
\nomenclature[AC]{$\Edges$}{set of edges (lines and transformers). 
\nomenclature[AD]{$\Disasters$}{set of disaster scenarios.}%
\nomenclature[AE]{$\Damages_s$}{set of edges that are inoperable during $s\in S$.}
\nomenclature[AF]{$\HardenedDamages_s$}{set of edges that are inoperable even though they are hardened during disaster $s\in S$.}
\nomenclature[AG]{$\lineCost_{ij}$}{cost to build a line between bus $i$ and $j$. 0 if line already exists.}
\nomenclature[AH]{$\switchCost_{ij}$}{cost to build a switch on a line between bus $i$ and $j$.}
\nomenclature[AI]{$\hardenCost_{ij}$}{cost to harden a line between bus $i$ and $j$.}
\nomenclature[AJ]{$\microgridCost_{i,k}$}{cost of generation capacity on phase $k$ at bus $i$.}
\nomenclature[AK]{$\facilityCost_{i}$}{cost to build a generation facility at node $i$.}
\nomenclature[AL]{$\capacity_{ijk}$}{line capacity between bus $i$ and bus $j$ on phase $k$.}
\nomenclature[AM]{$\phases_{ij}$}{set of phases for the line between bus $i$ and bus $j$.}
\nomenclature[AN]{$\phases_{i}$}{set of phases allowed to consume or inject at bus $i$.}
\nomenclature[AO]{$\phaseVariation_{ij}$}{parameter for controlling how much variation in flow between the phases is allowed.}
\nomenclature[AP]{$\demand_{i,k}$}{demand for power at bus  $i$ for phase $k$.}
\nomenclature[AQ]{$\availablecapacity_{i,k}$}{existing generation capacity on phase $k$ at node $i$.}
\nomenclature[AR]{$\microgridcapacity_{i,k}$}{maximum amount of generation capacity on phase $k$ that can be built at node $i$.}
\nomenclature[AS]{$\cycles$}{the set of sets of nodes that includes a cycle.}
\nomenclature[AT]{$\criticalLoad$}{fraction of critical load that must be served.}
\nomenclature[AU]{$\load$}{fraction of all load that must be served.}
\nomenclature[AV]{$\Load$}{set of buses whose load is critical.}

\nomenclature[BA]{$\textbf{Variables}$}{}%
\nomenclature[BB]{$\lineVariable_{ij}$}{determines if line $i,j$ is built.}%
\nomenclature[BC]{$\switchVariable_{ij}$}{determines if line $i,j$ has a switch.}
\nomenclature[BD]{$\hardenVariable_{ij}$}{determines if line $i,j$ is hardened.}
\nomenclature[BE]{$\microgridVariable_{i,k}$}{determines the capacity for generation on phase $k$ at node $i$.}
\nomenclature[BF]{$\facilityVariable_{i}$}{determines the generation capacity built at node $i$.}
\nomenclature[BG]{$\lineuseVariable_{ij}$}{determines if line $i,j$ is used during disaster $s$.}
\nomenclature[BH]{$\switchuseVariable_{ij}$}{determines if switch $i,j$ is used during disaster $s$.}
\nomenclature[BI]{$\linehardenVariable_{ij}$}{determines if line $i,j$ is hardened during disaster $s$.}
\nomenclature[BJ]{$\microgriduseVariable_{i,k}$}{determines the capacity for generation on phase $k$ at bus $i$ during disaster $s$.}
\nomenclature[BK]{$\facilityuseVariable_{i}$}{indicates if the generation capacity is used at node $i$ during disaster $s$.}
\nomenclature[BL]{$\generatorVariable_{i,k}$}{generation produced for bus $i$ on phase $k$ during disaster $s$.}
\nomenclature[BM]{$\loadVariable_{i,k}$}{load delivered at bus $i$ on phase $k$ during disaster $s$.}
\nomenclature[BN]{$\loadServeVariable_{i_j}$}{determines if the $j$th load at bus $i$ is served or not during disaster $s$.}
\nomenclature[BO]{$\flowVariable_{ij,k}$}{flow between bus $i$ and bus $j$ on phase $k$ during disaster $s$.}
\nomenclature[BP]{$\lineCycleVariable_{ij}$}{determines if at least one edge between $i$ and $j$ is used during disaster $s$.}
\nomenclature[BQ]{$\switchCycleVariable_{ij}$}{determines if at least one switch between $i$ and $j$ is used during disaster $s$.}
\nomenclature[BR]{$\lineuseVariable_{ij,0}$}{determines if there exists flow on line $i,j$ from $j$ to $i$, during disaster $s$.}
\nomenclature[BS]{$\lineuseVariable_{ij,1}$}{determines if there exists flow on line $i,j$ from $i$ to $j$, during disaster $s$.}

\printnomenclature

\paragraph{Distribution Grid Modeling} A distribution network is modeled as graph with nodes $\Nodes$ (buses) and edges $\Edges$ (power lines and transformers). In the physical system, each edge is composed of one, two, or three circuits or ``phases" and the electrical loads at the nodes are connected to and consume power from specific phases \cite{Garcia2000} (${\cal{P}}$). In many papers, multiple phases are approximated as a single phase with a single edge flow. However, under the damaged and stressed conditions considered in this work, the flows on the phases are often unbalanced, i.e. unequal, making it important to model all phases to accurately evaluate flow constraints on each phase. The phase flows are {\it not} directly controllable, but are related to nodal voltages and power injections by non-convex, physics-based equations \cite{Garcia2000}. Incorporation of these equations into the current formulation increases the complexity, however, the structure of distribution networks enables a simplification.

The design of protection systems for the vast majority of distribution circuits is based on the these circuits having a tree-like structure.  Therefore, although distribution grids are often designed to contain many possible loops, switches are used to ensure that these grids are operated in a tree or forest topology.  While, the switches introduce binary variables that increase the complexity of the ORDGDP, a linearized version of the electrical power flow equations (i.e. DC power flow) on the resulting trees is equivalent to a commodity flow model. We use a multi-commodity flow model that models each phase separately (Fig. \ref{fig:subproblem_feasible_set}).

The linearization of the power flow equations assumes uniform voltage magnitude at all nodes and ignores reactive power flows. In practice, we expect these are reasonable approximations because, prior to being upgraded, the distribution grid is already feasible with respect to voltage and reactive power flows.  By adding lines or distributed power sources, we put loads closer to generation thereby reducing voltage variability and reactive power flow and the potential for violating unmodeled constraints. In principle, it is possible to construct solutions where this is not the case, but the solutions to ORDGDP found by our algorithms has not resulted in these situations. However, this is an important area of future work, and we are developing methods to eliminate solutions that violate voltage or reactive power flow limits.

\paragraph{Damage Modeling} The ORDGDP is also defined by a set of scenarios, $\Disasters$. These scenarios are provided by a user or are drawn from a probabilistic damage model (the case here).  Each scenario is defined by the lines of the network that are damaged and are inoperable.

\paragraph{Design Options} We focus on four user-definable design options in distribution networks: 1) Hardening existing lines to lower the probability of damage, 2) Build new lines to add redundancy, 3) Build switches, to add operating flexibility, and 4) building distributed generation (sources of power).

\begin{figure}[!ht]
\tiny
\begin{alignat}{2}
\SubProblem = & \{\lineuseVariable,\switchuseVariable,\linehardenVariable,\microgriduseVariable,\facilityuseVariable : & & \nonumber \\
& -\lineuseVariable_{ij,0} \capacity_{ijk} \leq \flowVariable_{ij,k} \leq \lineuseVariable_{ij,1} \capacity_{ijk} & \forall ij\in \Edges, k \in \phases_{ij} & \label{eqn:line_constraint} \\
& \lineuseVariable_{ij,0} + \lineuseVariable_{ij,1} \leq \lineuseVariable_{ij} & \forall ij\in \Edges & \label{eqn:line_direction_constraint} \\
& (\switchuseVariable_{ij}-1) \capacity_{ijk} \leq \flowVariable_{ij,k} \leq (1-\switchuseVariable_{ij}) \capacity_{ijk} \quad & \forall ij\in \Edges, k \in \phases_{ij} & \label{eqn:switch_constraint} \\
& \frac{\displaystyle\sum_{k\in \phases_{ij}}f_{ij,k}}{\frac{|\phases_{ij}|}{(1-\beta_{ij})}} \leq \flowVariable_{ij,k^\prime} \leq \frac{\displaystyle\sum_{k\in \phases_{ij}}f_{ij,k}}{\frac{|\phases_{ij}|}{(1+\beta_{ij})}}\quad & \forall ij\in \Edges, k^\prime \in \phases_{ij} & \label{eqn:variation_constraint}\\
& \lineuseVariable_{ij} = \linehardenVariable_{ij} \leq \begin{cases} 0 & \text{ if } ij\in \HardenedDamages_s \\ 1 & \text{ else} \end{cases} & \forall ij\in \Damages_s & \label{eqn:damage}\\
& \loadVariable_{i,k} = \sum_{j = 0}^{n_i} \loadServeVariable_{i_j} \demand_{i_j,k} & \forall i\in \Nodes, k\in \phases_i & \label{eqn:load_constraint}\\
& 0 \leq \generatorVariable_{i,k} \leq \microgriduseVariable_{i,k} + \availablecapacity_{i,k} & \forall i\in \Nodes, k\in \phases_i & \label{eqn:generation_constraint}\\
& \generatorVariable_{i,k} - \loadVariable_{i,k} - \sum_{j\in \Nodes} \flowVariable_{ij, k} = 0 & \forall i\in \Nodes, k\in \phases_i & \label{eqn:balance_constraint}\\
& 0 \leq \microgriduseVariable_{i,k} \leq \microgridcapacity_{i,k}\facilityVariable_{i} & \forall i\in \Nodes, k\in \phases_i & \label{eqn:microgrid_constraint}\\
& \sum_{ij\in \Edges(C)} (\lineuseVariable_{ij} - \switchuseVariable_{ij}) \leq |V| - 1 & \forall C\in \cycles & \label{eqn:cycle_constraint}\\
& \switchuseVariable_{ij} \leq \lineuseVariable_{ij} & \forall ij\in \Edges & \label{eqn:switchable}\\
& \sum_{i \in \Load, k\in \phases_i} \loadVariable_{i,k} \geq \criticalLoad \sum_{i \in \Load, k\in \phases_i} \demand_{i,k} & & \label{eqn:critical_constraint}\\
& \sum_{i \in \Nodes \setminus \Load , k\in \phases_i} \loadVariable_{i,k} \geq \load \sum_{i \in \Nodes \setminus \Load , k\in \phases_i} \demand_{i,k} & & \label{eqn:loadserve_constraint}\\
& \lineuseVariable,\loadServeVariable, \switchuseVariable, \facilityuseVariable, \linehardenVariable \in \{0,1\}\} & & \label{eqn:discrete}
\end{alignat}

\caption{Set of feasible distribution networks}
\label{fig:subproblem_feasible_set}

\end{figure}

\paragraph{Optimization model}
Given a disaster $s\in S$, $\SubProblem$ in Fig. \ref{fig:subproblem_feasible_set} defines the set of feasible distribution networks.
The constraints of $\SubProblem$ involve a number of well-known constraints in the combinatorial optimization literature, including knapsacks, multi commodity flows, and tree constraints.
In this model, Eq. \ref{eqn:line_constraint} is a capacity constraint on phase flows.
When the line is not built the flow is forced to 0 by $\lineuseVariable$.
Eq. \ref{eqn:line_direction_constraint} forces all phases to flow in the same direction, an engineering constraint.
Eq. \ref{eqn:switch_constraint} states that the flow on a line is 0 when the switch is open.
Eq. \ref{eqn:variation_constraint} limits the fractional flow imbalance between the phases to smaller than $\phaseVariation_{ij}$.
Imbalance between phases cannot be extreme otherwise equipment may be damaged.
Here, we use $\phaseVariation_{ij}$ = 0.15 for transformers, and 1.0 otherwise.
Eq. \ref{eqn:damage} removes components in the damage set from the network by linking the two damage sets with the hardening variables.
Eq. \ref{eqn:load_constraint} requires that all or none of the load at a bus is served. Once again, this an engineering limitation of most networks.
Eq. \ref{eqn:generation_constraint} limits the distributed generation output by the generation capacity.
Eq. \ref{eqn:balance_constraint} ensures flow balance at the nodes for all phases.
Eq. \ref{eqn:microgrid_constraint} caps the generation capacity installed at the nodes.
Eq. \ref{eqn:cycle_constraint} eliminates network cycles, forcing a tree or forest topology.
Eq. \ref{eqn:switchable} states a switch is used only if the line exists.
Eq. \ref{eqn:critical_constraint} ensures a minimum fraction $\criticalLoad$ of critical load is served. Here, we generally require $\criticalLoad$ = 0.98.
Eq. \ref{eqn:loadserve_constraint} ensures that a minimum fraction of load is served. Here, $\load = 0.5$.  Eqs~\ref{eqn:critical_constraint} and \ref{eqn:loadserve_constraint} are the resilience criteria that must be met by  $\SubProblem$ and are similar to the $n-k-\epsilon$ criteria of \cite{Chen2014a}.
Eq. \ref{eqn:discrete} states which variables are discrete.

One of the more difficult constraints in this formulation is Eq. \ref{eqn:cycle_constraint} due to possible combinatorics. There are a number of ways to implement cycle constraints and we use the formulation in Fig. \ref{fig:cycle_trick}.


\begin{figure}[!ht]
\begin{alignat}{2}
& \sum_{ij\in \Edges(C)} (\lineCycleVariable_{ij} - \switchCycleVariable_{ij}) \leq |V| - 1 & \forall C\in \cycles & \label{eqn:no_cycle_1}\\
& \lineuseVariable_{ij} \leq \lineCycleVariable_{ij} & \forall ij \in \Edges & \label{eqn:line_implication_constraint} \\
& 3 - \lineuseVariable_{ij} - \switchCycleVariable_{ij} \geq \switchuseVariable_{ij} \geq \lineuseVariable_{ij} + \switchCycleVariable_{ij} - 1\quad & \forall ij \in \Edges & \label{eqn:switch_implication_constraint}
\end{alignat}
\caption{Cycle constraints}
\label{fig:cycle_trick}
\end{figure}

\noindent While the multi-graph structure introduces a large number of cycles, there are a relatively small number of cycles when the multi-edges are reduced to one edge. Thus, we introduce binary variables (linear number) for the edges of the corresponding single-edge graph and enumerate the possible cycles in that graph (Eq. \ref{eqn:no_cycle_1}).  Then, Eqs \ref{eqn:line_implication_constraint} and \ref{eqn:switch_implication_constraint} are used to pass information between artificial cycle variables and the actual line and switch variables.

For each $s \in \Disasters$, $\SubProblem$ determines the set of feasible distribution networks for each $s$. There are some redundant variables in this formulation that improves the separability of the problem. The ORDGDP is the minimum cost design that falls in the intersection of all the $\SubProblem$  (Fig. \ref{fig:optimization_model}).

\begin{figure}[!ht]
\tiny
\begin{alignat}{2}
\min\quad &\sum_{ij\in \Edges} \lineCost_{ij} \lineVariable_{ij} + \sum_{ij \in \Edges} \switchCost_{ij} \switchVariable_{ij} + \sum_{ij \in \Edges} \hardenCost_{ij} \hardenVariable_{ij}\nonumber \\
& + \sum_{i \in \Nodes} \facilityCost_{i} \facilityVariable_{i} + \sum_{i \in \Nodes, k\in \phases_i} \microgridCost_{i,k} \microgridVariable_{i,k} \label{eqn:obj}\\
\text{s.t.}\quad & \lineuseVariable_{ij} \leq \lineVariable_{ij} & \forall ij\in \Edges, s\in \Disasters & \label{eqn:assign:1}\\
& \switchuseVariable_{ij} \leq \switchVariable_{ij} & \forall ij\in \Edges, s\in \Disasters & \label{eqn:assign:2} \\
& \linehardenVariable_{ij} \leq \hardenVariable_{ij} & \forall ij\in \Edges, s\in \Disasters & \label{eqn:assign:3} \\
& \microgriduseVariable_{i,k} \leq \microgridVariable_{i,k} & \forall i\in \Nodes, k\in \phases_i, s\in \Disasters & \label{eqn:assign:4} \\
& \facilityuseVariable_{i} \leq \facilityVariable_{i} & \forall i\in \Nodes, s\in \Disasters & \label{eqn:assign:5} \\
& \microgridVariable_{i,k} \leq M_{i,k} \facilityVariable_{i} & \forall i\in \Nodes, k\in \phases_i & \label{eqn:assign:6} \\
& (\lineuseVariable,\switchuseVariable,\linehardenVariable,\microgriduseVariable,\facilityuseVariable) \in \SubProblem & \forall s\in \Disasters \label{eqn:feasible_network} \\
& \lineVariable,\switchVariable, \hardenVariable, \facilityVariable \in \{0,1\} & & \label{eqn:discrete_2}
\end{alignat}
\caption{Optimal Resilient Distribution Grid Design}
\label{fig:optimization_model}
\end{figure}

\noindent Eq. \ref{eqn:obj} minimizes the cost of building lines and switches, hardening lines, and building facilities and generation.
For notational simplicity, existing lines, switches, and generation are included as variables in the objective with 0 cost, however in practice these enter the formulation as constants.
Eqs. \ref{eqn:assign:1} through \ref{eqn:assign:6} tie the first stage (construction) decisions with second stage variables ($\SubProblem$). Eq. \ref{eqn:feasible_network} states that the mixed-integer vector $(\lineuseVariable,\switchuseVariable,\linehardenVariable,\microgriduseVariable,\facilityuseVariable)$ constitutes a feasible distribution network for scenario $s$.

\paragraph{Chance Constraints}
For some networks, a very small number of scenarios in $\Disasters$ may drive the total cost in Eq. \ref{eqn:obj}. In real-world applications, the designer of the network may lower the total investment cost by accepting some risk of not always satisfying the resiliency criteria.  In these situations, we can relax Eqs. \ref{eqn:critical_constraint} and \ref{eqn:loadserve_constraint} to a set of chance constraints:
\begin{equation}
\tiny
P\left(\begin{array}{ll}
\sum_{i \in \Load, k\in \phases_i} \loadVariable_{i,k} \geq \criticalLoad \sum_{i \in \Load, k\in \phases_i} \demand_{i,k} & \forall s\in \Disasters\\
\sum_{i \in \Nodes, k\in \phases_i} \loadVariable_{i,k} \geq \load \sum_{i \in \Nodes, k\in \phases_i} \demand_{i,k} & \forall s\in \Disasters
\end{array}\right) \geq 1-\epsilon
\end{equation}

When assuming the scenarios follow a uniform distribution, this is equivalent to stating that these constraints are violated in $\epsilon|S|$ of the scenarios. Thus, we can restate these constraints as:
\begin{equation}\label{eq:CC}
\tiny
\begin{array}{ll}
\sum_{i \in \Load, k\in \phases_i} \loadVariable_{i,k} \geq \criticalLoad \sum_{i \in \Load, k\in \phases_i} \demand_{i,k} (1-\ChanceVariable_s) & \forall s\in \Disasters\\
\sum_{i \in \Nodes, k\in \phases_i} \loadVariable_{i,k} \geq \load \sum_{i \in \Nodes, k\in \phases_i} \demand_{i,k} (1-\ChanceVariable_s) & \forall s\in \Disasters\\
\sum_{s\in S} \ChanceVariable_s \leq \epsilon |S|
\end{array}
\end{equation}

\section{Algorithms}

In this section we discuss the algorithms we developed for solving the ORDGDP. ORDGDP is a two-stage mixed integer programming (MIP) problem with a block diagonal structure that includes coupling variables between the blocks.  We developed an exact algorithm that is vastly more efficient than a commercial state-of-the-art MIP solver.  We then used the exact algorithm to develop a hybrid with variable neighborhood search that is competitive with the exact solver and is better than a heuristic used by the industry.

\paragraph{Scenario-Based Decomposition (SBD)}

Decomposition is often used for solving two-stage stochastic MIPs, and it can be applied to ORDGDP after the following key observation:

\begin{obs}
The second stage variables do not appear in the objective function, therefore any optimal first stage solution based on a sub set of the second stages and is feasible for all second stage sub problems is an optimal solution.
\end{obs}

Based on this observation we can apply SBD to solve the ORDGDP. At high level, this algorithm solves problems with iteratively larger sets of scenarios until a solution is obtained that is feasible for all scenarios. The algorithm takes as input the set of disasters (scenarios) and an initial scenario to consider, $S^\prime$.  Line 2 solves ORDGDP on $S^\prime$,
where $P(S^\prime)$ and $\sigma^*$ are used to denote the problem and solution respectively.  Line  3 then evaluates   $\sigma^*$ on the remaining scenarios in $S \setminus S^\prime$.
The function $l:P^\prime(s,\sigma^*) \rightarrow \mathbb{R}_+$, is an infeasibility measure that is 0 if the problem is feasible, positive otherwise. This is implemented by maximizing the reliability constraints, i.e. total and critical demand satisfied. It measures the gap between the delivered and the required demand (the right hand side of the Eqs. \ref{eqn:critical_constraint} and \ref{eqn:loadserve_constraint}). This function prices the current solution over $s\in S\setminus S^\prime$. If all prices are 0, then the algorithm terminates with solution $\sigma^*$ (lines 4-5).  Otherwise, the algorithm adds the scenario with the worst infeasabilty measure to $S^\prime$ (line 7).


This scenario-based decomposition shares a number of key features with Benders decomposition \cite{VW:10}.  Both Benders and SBD solve the subproblems based on iterative solves to the first-stage problem.  Benders will typically add a single constraint in the form of cut that represents a facet of the subproblem, whereas in SBD we add the entire polyhedron associated with the sub-problem.
Benders is a successful approach on similar expansion, commodity flow, and interdiction problems, but is unsuccessful here due to the non-convex (discrete) nature of the second stage. Hence our generalization to SBD.

\begin{algorithm}[!ht]
\small
\SetKwInOut{Input}{input}
\Input{A set of disasters $S$ and let $S^\prime = S_0$\;}
\While{$S\setminus S^\prime \neq \emptyset$}{
	$\sigma^* \leftarrow$ Solve $P(S^\prime)$\;
	$I \leftarrow \left<s_1, s_2 \ldots s_{|S\setminus S^\prime|} \right> s\in S\setminus S^\prime :$ $l(P^\prime(s_i,\sigma^*)) \ge l(P^\prime(s_{i+1},\sigma^*))$\;
	\eIf{$l(P^\prime(I(0),\sigma^*)) \leq 0$}{
		\Return{$\sigma^*$}\;
	}{
		$S^\prime \leftarrow S^\prime \cup  I(0)$\;
	}
\Return{$\sigma^*$}
}
\caption{Scenario Based Decomposition}
\label{benders}
\end{algorithm}

\paragraph{Greedy Algorithm}
A computationally efficient way of generating feasible solutions to the ORDGDP relaxes the coupling first stage variables and solves each scenario $s \in \Disasters$ individually.
The solutions are combined by taking the maximum of each construction variable (${\cal{X}} = \lineVariable \cup \switchVariable \cup \hardenVariable \cup \microgridVariable \cup \facilityVariable$) over all scenarios (Algorithm 2).
The switch construction cost is determined by switches that are needed to reduce the network into a tree for every scenario (line 4). Although the Greedy Algorithm is simple and fast, it rarely results in an optimal investment decision. However, it is representative of the types of heuristics used by the industry: see Reference \cite{Munoz2014} for a survey.

\begin{algorithm}[!ht]
\small
\SetKwInOut{Input}{input}
\Input{A set of disasters $S$\;}
\For{$s\in S$}{
	$\sigma^s \leftarrow Solve(P^\prime(s))$\;
}
$\sigma^*(x) = \max\{\sigma^s(x) | \forall s\in S\},\; \forall x \in {\cal{X}}$\;
Update $\sigma^*(x_i)$ with switches to preserve feasibility\;
\Return{$\sigma^*$}
\caption{Greedy}
\label{greedy}
\end{algorithm}

\paragraph{Variable Neighborhood Search}
To overcome the limitations of greedy heuristics like Algorithm 2, we developed an approach based on Variable Neighborhood Decomposition (VNS) Search \cite{LHMU:10}.
The algorithm fixes a subset of first stage variables to their current value and searching the remaining variables for a better solution. If all the first stage variables are fixed, the problem decomposes into $|S|$ separate problems that are easily solved and provide heuristic justification for focusing on first stage variables. More formally, $P(\sigma,J)$ denotes the problem with first stage variables, $J \in  {\cal{X}}$,fixed to $\sigma$, i.e. $x_j = \sigma(x_j)$, and $P^{LP}$ is the LP relaxation of problem $P$.

Algorithm 3 describes the VNS procedure. Line  1 computes the solution to the LP relaxation of the ORDGDP, $(\sigma^{LP})$. Line 4 counts the number of variable assignments that are different between the solution to LP relaxation $(\sigma^{LP})$ and the best known solution $\sigma^*$ ($\sigma(x)$ denotes the variable assignment of $x$ in solution $\sigma$). Line 5 orders the variables of $\cal{X}$ by the difference between their assignments in $\sigma^*$ and $\sigma^{LP}$.  Heuristically, those variables whose assignments are furthest from their LP assignment represent good opportunities to improve $\sigma^*$. The algorithm updates the rate at which the neighborhood size is increased ($step$) based on whether or not the algorithm is in a restart situation (lines 8 and 11).  If the algorithm is in a restart, the ordering of the variables are also randomized (line 9).  Line 13 computes the best solution in the neighborhood of $\sigma$ where the first $k$ elements of $J$ are fixed.  If the resulting solution is better, then the algorithm proceeds with a new $\sigma^*$ (lines 15-18)--$f$ is used as shorthand for Eq.~\ref{eqn:obj}. Otherwise, the size of the neighborhood is increased (lines 20-23). The iterations terminate when the maximum number of restarts is reached (line 2), the maximum number of neighborhood resizings is reached (line 12), or a time limit is reached. In this paper, $\textsc{maxRestarts} = 10$, $\textsc{maxIterations} = 4$, $\textsc{maxTime} = 48$ CPU hours, and $d = 2$.

\begin{algorithm}[!ht]
\small
\SetKwInOut{Input}{input}
\Input{$\sigma^\prime$, $\textsc{maxTime}$, $\textsc{maxRestarts}$ and $\textsc{maxIterations}$\;}
Let $\sigma^{LP}\leftarrow Solve(P^{LP})$, $\sigma^* \leftarrow \sigma^\prime$, $restart \leftarrow false$;\\
\While{$t < \textsc{maxTime}$ and $i < \textsc{maxRestarts}$}{
    $j \leftarrow 0$\;
	$n \leftarrow | x \in \cal{X} : |\sigma^*($$x$$) - \sigma^{LP}($$x$$)| \ne 0|$\;
	$J \leftarrow \left<\pi_1, \pi_2 \ldots \pi_{|J|} \right> \in  \cal{X} :$ $|\sigma^*(\pi_i) - \sigma^{LP}(\pi_i)| \le |\sigma^*(\pi_{i+1}) - \sigma^{LP}(\pi_{i+1})|$\;
	\eIf{$restart$}{	
		$i \leftarrow i + 1$\;
		$step \leftarrow \frac{4n}{d}$, $k = |{\cal{X}}| - step$\;
		\textit{shuffle}($J$)
	}{
		$step \leftarrow \frac{n}{d}$, $k = |{\cal{X}}| - step$\;
	}
	\While{$t < \textsc{maxTime}$ and $j \le \textsc{maxIterations}$}{
		$\sigma^\prime \leftarrow Solve(P(\sigma^*, J(1,\ldots,k))$\;
		\eIf{$f(\sigma^\prime) < f(\sigma^*)$}{
			$\sigma^* \leftarrow \sigma^\prime$\;
			$i \leftarrow 0$\;
			$restart \leftarrow false$\;
			$j \leftarrow \textsc{maxIterations}$\;
		}{
            $j \leftarrow j+1$\;
			$k = k - \frac{step}{2}$\;
			\If{$j > \textsc{maxIterations}$}{
				$restart \leftarrow true$\;
			}
		}
	}	
\Return{$\sigma^*$}
}
\caption{Variable Neighborhood Search}
\label{VNS}
\end{algorithm}


\paragraph{Scenario-based Variable Neighborhood Decomposition Search (SBVNDS)}
Given that we have a powerful exact method in Algorithm 1 as well as a VNS in Algorithm 3, the natural algorithm hybridizes these approaches to get Algorithm, SBVNDS.
The algorithm proceeds exactly the same same as Algorithm 1, except that the exact solver for $Solve(P(S^\prime))$ is replaced by VNS in line 2.


\section{Empirical Results}

\begin{table*}
\begin{tiny}
\begin{tabular}{cc}

\begin{tabular}{cc|c|c|c|c|c|c|}
& \multicolumn{7}{c}{Urban, Hardened lines are not damageable (a) } \\
\cline{2-8}
& \multicolumn{2}{|c|}{CPLEX} & \multicolumn{1}{c|}{Greedy} & \multicolumn{2}{c|}{SBD} & \multicolumn{2}{c|}{SBVNDS} \\
\cline{2-8} \cline{2-8}
& \multicolumn{1}{|c|}{CPU} & OBJ & OBJ & CPU & OBJ & CPU & OBJ \\
\cline{2-8} \cline{2-8}
$10\%$& \multicolumn{1}{|c|}{19984.7} & 322.9 & 1044.5 & 465.8 & 322.9 & 289.9 & 353.7 \\
\cline{2-8}
$25\%$& \multicolumn{1}{|c|}{166352} & 635.4 & 1643.5 & 8028.3 & 635.4 & 811.4 & 635.4 \\
\cline{2-8}
$50\%$& \multicolumn{1}{|c|}{TO} & X & 2021.2 & 2840.7 & 647.7 & 791.3 & 647.7 \\
\cline{2-8}
$75\%$& \multicolumn{1}{|c|}{TO} & X & 1874.2 & 991.1 & 652.1 & 692.5 & 652.1 \\
\cline{2-8}
$100\%$& \multicolumn{1}{|c|}{TO} & X & 1934.4 & 712.7 & 654.1 & 662.5 & 654.1 \\
\cline{2-8}
\end{tabular}

&

\begin{tabular}{cc|c|c|c|c|c|c|}
& \multicolumn{7}{c}{Rural, Hardened lines are not damageable (b)} \\
\cline{2-8}
& \multicolumn{2}{|c|}{CPLEX} & \multicolumn{1}{c|}{Greedy} & \multicolumn{2}{c|}{SBD} & \multicolumn{2}{c|}{SBVNDS} \\
\cline{2-8} \cline{2-8}
& \multicolumn{1}{|c|}{CPU} & OBJ & OBJ & CPU & OBJ & CPU & OBJ \\
\cline{2-8} \cline{2-8}
$10\%$& \multicolumn{1}{|c|}{33083.5} & 2337.0 & 3274.8 & 1837.9 & 2337.0 & 503.3 & 2337.0 \\
\cline{2-8}
$25\%$& \multicolumn{1}{|c|}{32170.8} & 2390.3 & 3427.6 & 571.0 & 2390.3 & 457.8 & 2390.3 \\
\cline{2-8}
$50\%$& \multicolumn{1}{|c|}{20840.3} & 2397.6 & 3449.9 & 471.2 & 2397.6 & 421.2 & 2397.6\\
\cline{2-8}
$75\%$& \multicolumn{1}{|c|}{15556.1} & 2400.4 & 3452.7 & 337.5 & 2400.4 & 299.8 & 2400.4 \\
\cline{2-8}
$100\%$& \multicolumn{1}{|c|}{17225.9} & 2400.6 & 2780.6 & 385.8 & 2400.6 & 346.9 & 2400.6 \\
\cline{2-8}
\end{tabular} \\
\\

\begin{tabular}{cc|c|c|c|c|c|c|}
& \multicolumn{7}{c}{Urban, Hardened lines are damaged at a $\frac{1}{100}$ rate (c) } \\
\cline{2-8}
& \multicolumn{2}{|c|}{CPLEX} & \multicolumn{1}{c|}{Greedy} & \multicolumn{2}{c|}{SBD} & \multicolumn{2}{c|}{SBVNDS} \\
\cline{2-8} \cline{2-8}
& \multicolumn{1}{|c|}{CPU} & OBJ & OBJ & CPU & OBJ & CPU & OBJ \\
\cline{2-8} \cline{2-8}
$10\%$& \multicolumn{1}{|c|}{159166} & 445.8 & 1061.7 & 2232.9 & 445.8 & 2721.3 & 476.5\\
\cline{2-8}
$25\%$& \multicolumn{1}{|c|}{TO} & X & 1441.9 & 14299.2 & 662.9 & 2994.7 & 701.5\\
\cline{2-8}
$50\%$& \multicolumn{1}{|c|}{TO} & X & 1571.2 & 2848.7 & 646.0 & 1917.7 & 760.2\\
\cline{2-8}
$75\%$& \multicolumn{1}{|c|}{TO} & X & 1787.3 & 16040.6 & 687.6 & 1481.4 & 687.6\\
\cline{2-8}
$100\%$& \multicolumn{1}{|c|}{TO} & X & 2744.8 & 24270.3 & 1320.5 & 2157.5 & 1330.5 \\
\cline{2-8}
\end{tabular}

&

\begin{tabular}{cc|c|c|c|c|c|c|}
& \multicolumn{7}{c}{Rural,  Hardened lines are damaged at a $\frac{1}{100}$ rate (d) } \\
\cline{2-8}
& \multicolumn{2}{|c|}{CPLEX} & \multicolumn{1}{c|}{Greedy} & \multicolumn{2}{c|}{SBD} & \multicolumn{2}{c|}{SBVNDS} \\
\cline{2-8} \cline{2-8}
& \multicolumn{1}{|c|}{CPU} & OBJ & OBJ & CPU & OBJ & CPU & OBJ \\
\cline{2-8} \cline{2-8}
$10\%$& \multicolumn{1}{|c|}{77947.9} & 2363.0 & 3375.4 & 759.0 & 2363.0 & 576.9 & 2363.0\\
\cline{2-8}
$25\%$& \multicolumn{1}{|c|}{TO} & X & 8238.6 & TO & X & 919.4 & 6744.3\\
\cline{2-8}
$50\%$& \multicolumn{1}{|c|}{TO} & X & 12336.0 & TO & 9288.9 & 4361.8 & 7121.0\\
\cline{2-8}
$75\%$& \multicolumn{1}{|c|}{TO} & X & 23099.5 & TO & X & 23142.6 & 11500.0\\
\cline{2-8}
$100\%$& \multicolumn{1}{|c|}{TO} & X & 16600.7 & TO & X & 5879.5 & 9797.3\\
\cline{2-8}
\end{tabular} \\

\\

\begin{tabular}{cc|c|c|c|c|c|c|}
& \multicolumn{7}{c}{Urban, Hardened lines are damaged at a $\frac{1}{10}$ rate (e)} \\
\cline{2-8}
& \multicolumn{2}{|c|}{CPLEX} & \multicolumn{1}{c|}{Greedy} & \multicolumn{2}{c|}{SBD} & \multicolumn{2}{c|}{SBVNDS} \\
\cline{2-8} \cline{2-8}
& \multicolumn{1}{|c|}{CPU} & OBJ & OBJ & CPU & OBJ & CPU & OBJ \\
\cline{2-8} \cline{2-8}
$10\%$& \multicolumn{1}{|c|}{TO} & X & 859.1 & 5265.1 & 460.8 & 2505.7 & 594.1 \\
\cline{2-8}
$25\%$& \multicolumn{1}{|c|}{TO} & X & 1742.2 & 12530.3 & 961.2 & 2843.2 & 961.2 \\
\cline{2-8}
$50\%$& \multicolumn{1}{|c|}{TO} & X & 3133.8 & 34822.7 & 1417.2 & 3363.5 & 1555.2 \\
\cline{2-8}
$75\%$& \multicolumn{1}{|c|}{TO} & X & 3472.0 & TO & X & 7486.5 & 1894.2\\
\cline{2-8}
$100\%$& \multicolumn{1}{|c|}{TO} & X & 10479.1 & TO & X & 32289.8 & 7959.4\\
\cline{2-8}
\end{tabular}

&

\begin{tabular}{cc|c|c|c|c|c|c|}
& \multicolumn{7}{c}{Rural, Hardened lines are damaged at a $\frac{1}{10}$ rate (f)} \\
\cline{2-8}
& \multicolumn{2}{|c|}{CPLEX} & \multicolumn{1}{c|}{Greedy} & \multicolumn{2}{c|}{SBD} & \multicolumn{2}{c|}{SBVNDS} \\
\cline{2-8} \cline{2-8}
& \multicolumn{1}{|c|}{CPU} & OBJ & OBJ & CPU & OBJ & CPU & OBJ \\
\cline{2-8} \cline{2-8}
$10\%$& \multicolumn{1}{|c|}{TO} & X & 7503.3 & 141718.0 & 4325.9 & 7756.8 & 4424.8\\
\cline{2-8}
$25\%$& \multicolumn{1}{|c|}{TO} & X & 18021.3 & TO & X & 21993.5 & 7371.9\\
\cline{2-8}
$50\%$& \multicolumn{1}{|c|}{TO} & X & 28865.0 & TO & 12017.7 & 74729.0 & 12031.2\\
\cline{2-8}
$75\%$& \multicolumn{1}{|c|}{TO} & X & 31887.0 & TO & 13522.2 & 107165.0 & 13500.8\\
\cline{2-8}
$100\%$& \multicolumn{1}{|c|}{TO} & X & 32901.9 & TO & 16794.4 & 114354.0 & 16778.2\\
\cline{2-8}
\end{tabular}

\end{tabular}

\end{tiny}

\caption{These tables compare the performance of the algorithms when hardened lines cannot be damaged (a, b),
are damaged at $\frac{1}{100}$ the rate of unhardened lines (c, d), and damaged at $\frac{1}{10}$ the rate of unhardened lines (e, f). The columns denoted by CPU and OBJ refer to CPU time and objective value, respectively. We omit the CPU time of Greedy as it is always less than 60 CPU seconds. The rows refer to the probability a 1 mile segment of a line is damaged. }
\label{table:tab1}
\end{table*}

The algorithms were implemented using the CPLEX C++ API with Concert technology as a 32 threaded application on Intel XEON 2.29 GHz processors.  Since these are planning problems, in principle, practitioners could utilize days of CPU time to produce a plan. However, in order to produce a wide range of results, we limited the algorithms to 48 hours of CPU time. Our problems are based on a modified version of the IEEE 34 bus systems \cite{Kersting1991} (see Fig.~\ref{fig:1}) that are representative of medium sized distribution systems. \footnote{Due to space constraints, details are omitted. We will purchase extra pages to provide more details}



\begin{figure}[!ht]
   \centering
   \subfigure[Urban]{\includegraphics[width=0.23\textwidth]{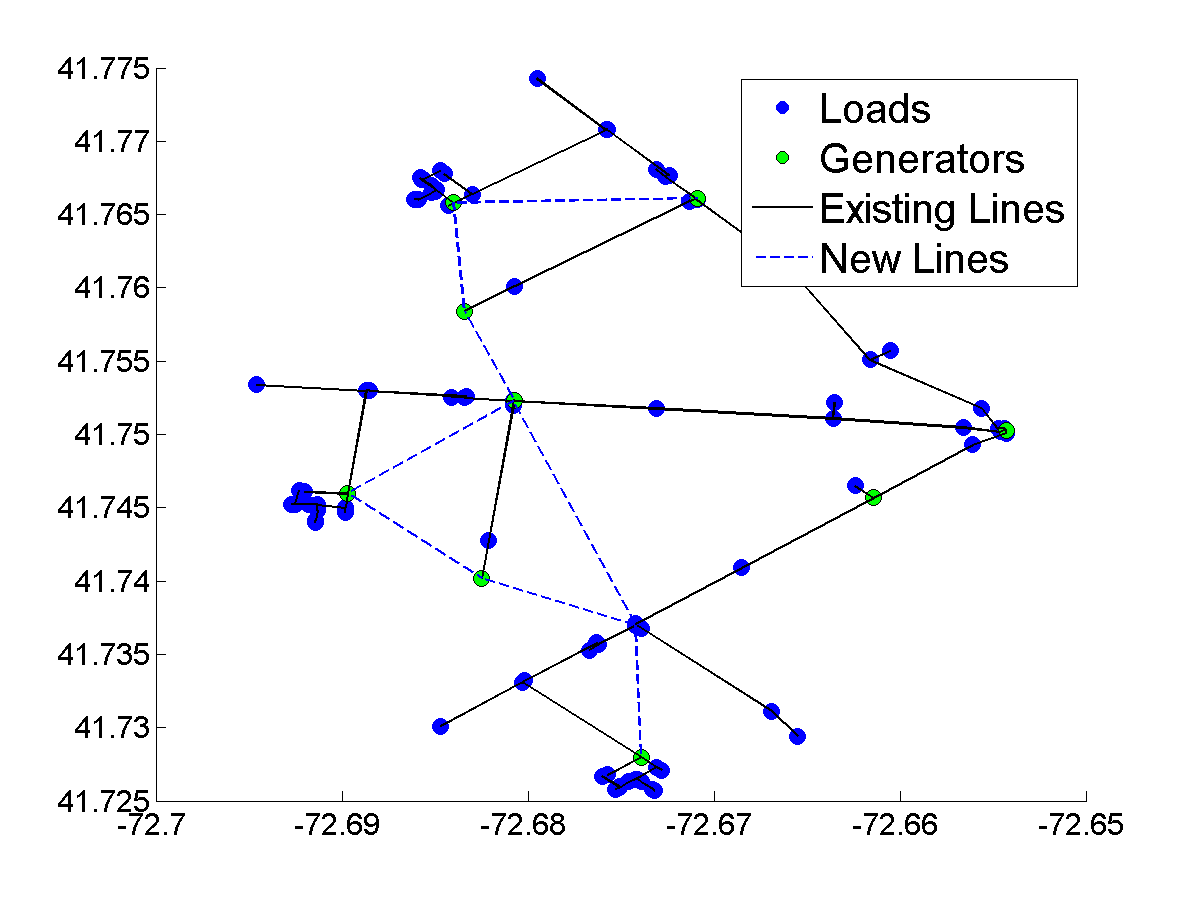}}\label{fig:1:a}
   \subfigure[Rural]{\includegraphics[width=0.23\textwidth]{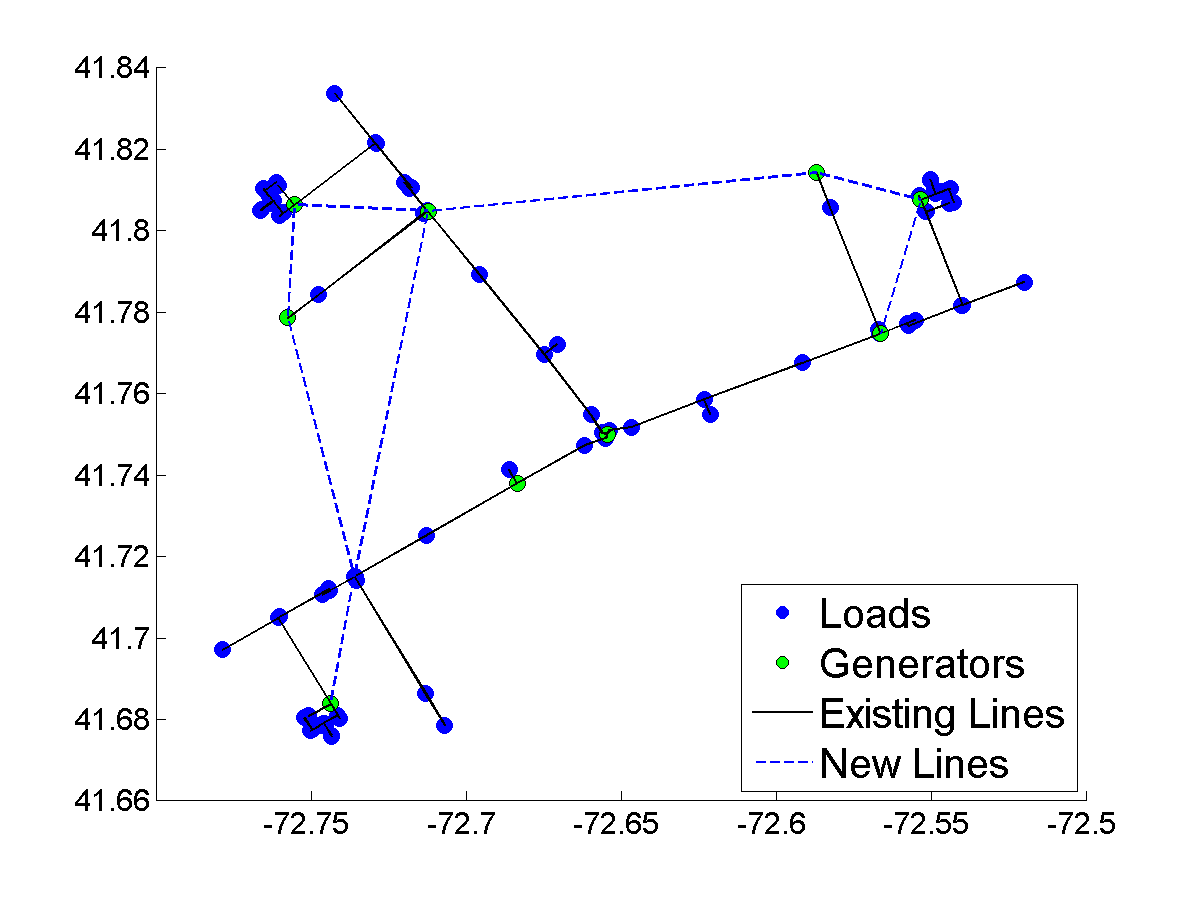}}\label{fig:1:b}

   \caption{We generated two variations of the IEEE 34 bus problem.  Each problem contains three copies of the IEEE 34 system to mimic situations where there are three normally independent distribution circuits that could support each other during extreme events.  These problems include 100 scenarios, 109 nodes, 118 possible generators, 204 loads, and 148 edges, resulting in problems with $> 90k$ binary variables. The difference between rural (a) and urban (b) is the distances between nodes (expansion costs). }
   \label{fig:1}
\end{figure}

Scenarios for this paper are based on damage caused by ice storms, whose intensity tends to be homogeneous on the scale of distribution systems \cite{Sa_Thesis}.  Intensities are modeled as damage rates per mile on power poles and are transformed into the probability a power line segment of one mile length is damaged (a pole has failed). Empirically, we find that 100 randomly created scenarios is sufficient to capture the salient features of the distribution. Each scenario contains two sets of line failures, one for hardened lines ($\HardenedDamages_s$) and a second for lines that are not hardened ($\Damages_s$).

Table \ref{table:tab1} provide results when hardened lines are not damaged or are damaged at rates of  $\frac{1}{100}$ or $\frac{1}{10}$ of the unhardened rate. There are a
number of important observations in these tables.  First, CPLEX by itself is computationally uncompetitive. Only when the hardened lines are not damaged, CPLEX completes within the time limit, but these problems are ``easier'' because hardened lines are robust and relatively inexpensive, enabling CPLEX to eliminate many solutions.  The objective for Greedy is always worse than optimal.  The exact method SBD is much more computationally efficient than CPLEX and is able to solve many more problems to optimality indicating that CPLEX is unable to recognize the scenario structure in the problems.  However, SBD is sensitive to which scenarios are included (function $l$), and if poor choices are made, it begins to resemble CPLEX. However, the meta-heuristic SBVNDS is able to overcome these limitations.  It is much faster than SBD, and almost always achieves the optimal solution.  This indicates that heuristic methods based on combining powerful techniques like VNS with strong exact algorithms are very good on this type of 2-stage mixed integer programming problems.

\paragraph{Critical load constraint}
Figures \ref{fig:12} and \ref{fig:13} show some results for rural and urban problems when the required fraction of critical load served is varied. In general, peaks in CPU time correspond to discrete jumps in the amount of load served as $\lambda$ increases.


\begin{figure}[!ht]
   \centering
   \subfigure[CPU time]{\includegraphics[width=0.23\textwidth]{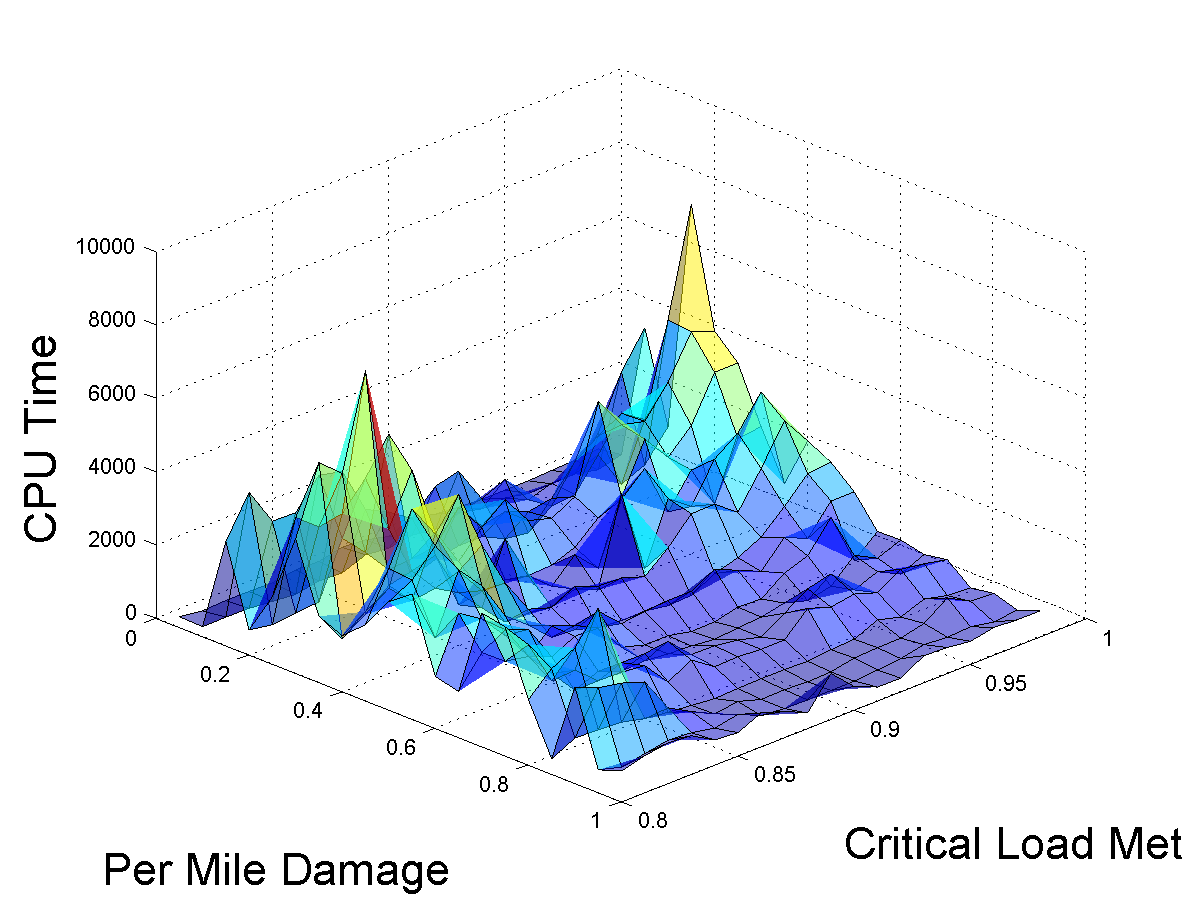}}\label{fig:12:a}
   \subfigure[Objective value]{\includegraphics[width=0.23\textwidth]{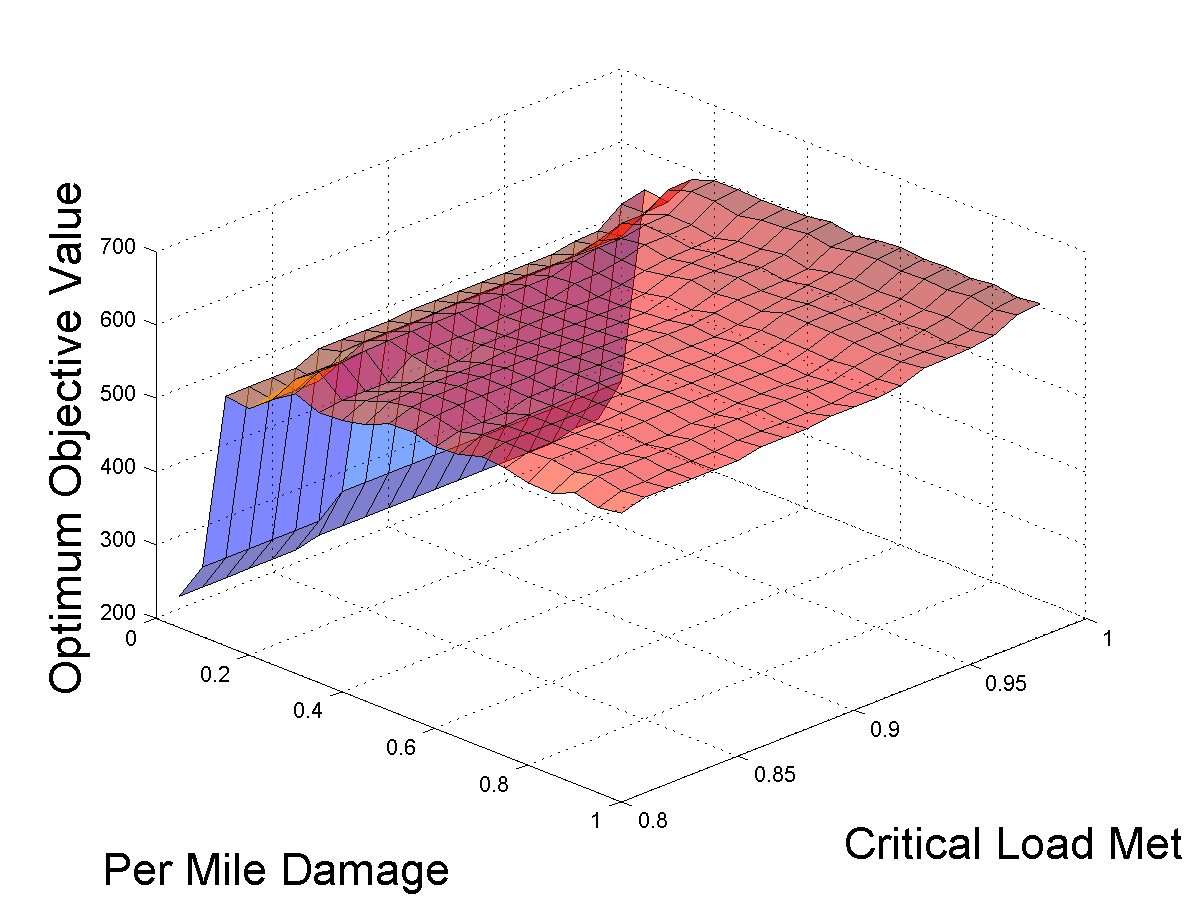}}\label{fig:12:b}
   \caption{Sensitivity of the CPU time and objective value to changes in $\lambda$ on the Urban problem for SBD when hardened lines are not damageable. Due to short distances, the solution favors hardening many lines. The required hardening is relatively insensitive to the amount of damage and $\lambda$. However, there are spikes in problem difficulty at transitions in $\lambda$ that require additional load service.}
   \label{fig:12}
\end{figure}


\begin{figure}[!ht]
   \centering
   \subfigure[CPU time]{\includegraphics[width=0.23\textwidth]{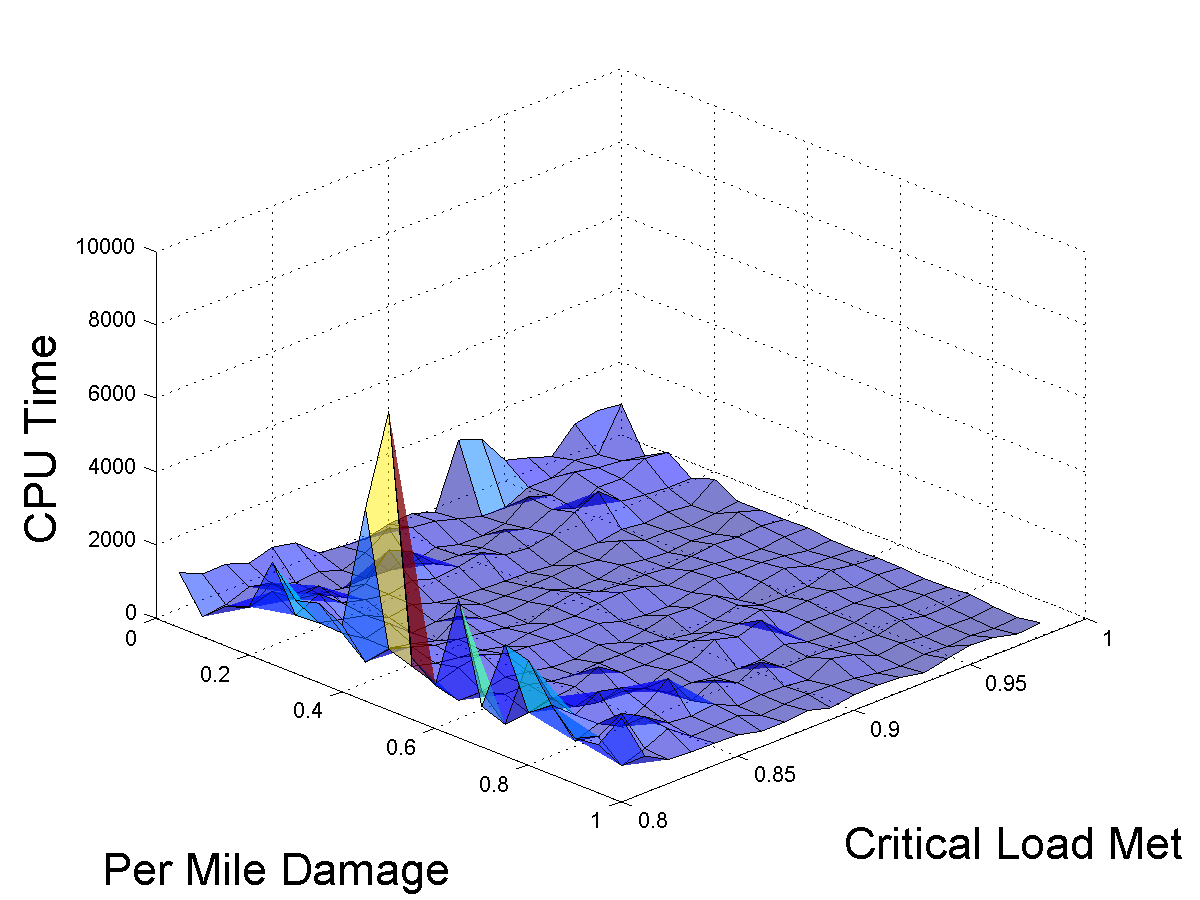}}\label{fig:13:a}
   \subfigure[Objective value]{\includegraphics[width=0.23\textwidth]{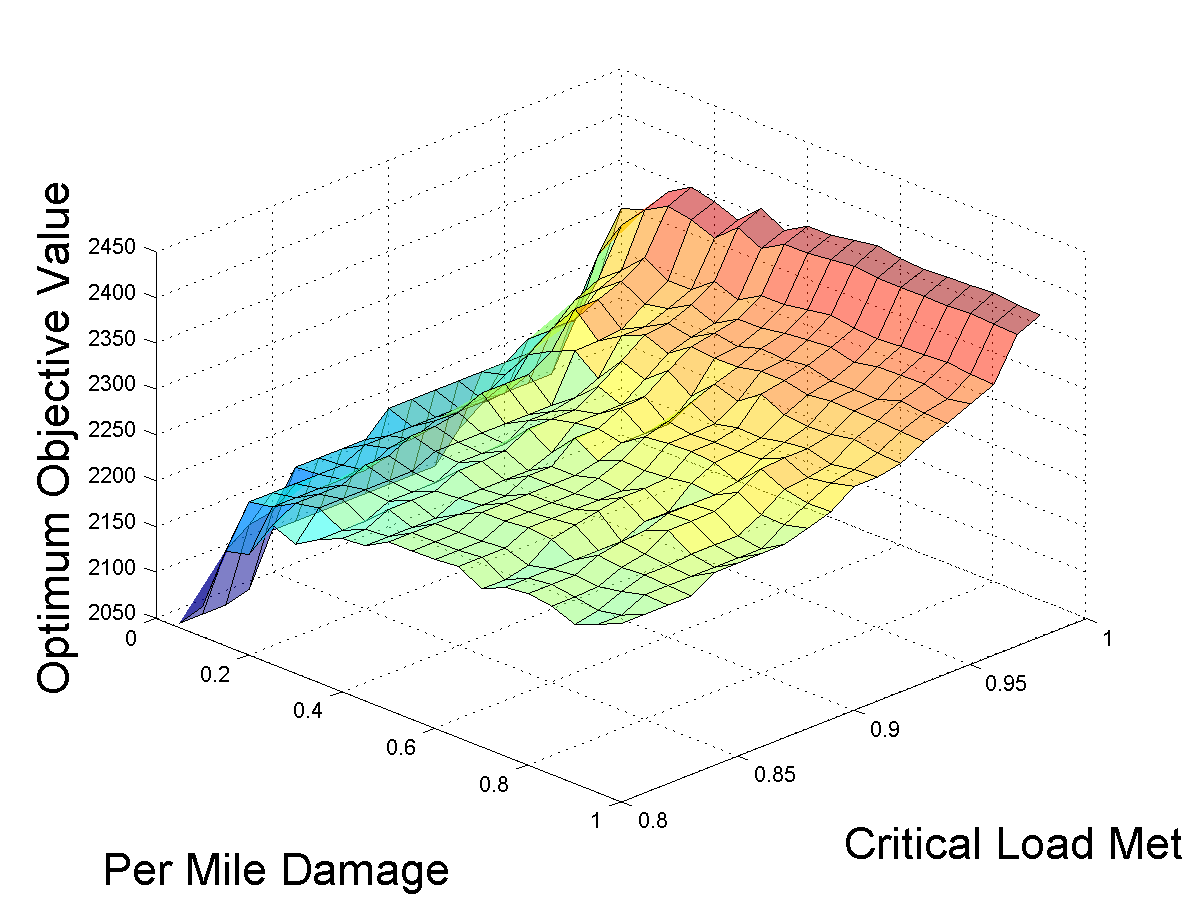}}\label{fig:13:b}
   \caption{Sensitivity of the CPU time and objective value to changes in $\lambda$ for SBD on the Rural problem when hardened lines are not damageable. Because of long distances, the solution favors adding generation and is sensitive to the amount of damage and $\lambda$. }

   \label{fig:13}
\end{figure}


\paragraph{Chance constraints}
Fig.~\ref{fig:14} shows results when the resiliency criteria are relaxed to the chance constraints in Eq.~\ref{eq:CC} and $\epsilon$ is varied. Interestingly, CPU time is not impacted too greatly by damage rates. Also, the solution is relatively insensitive to the choice of $\epsilon$ as damage rates increase, indicating that an ``easier'' problem with small $\epsilon$ could be used to approximate a solution to the harder problems.


\begin{figure}[!ht]
   \centering
   \subfigure[CPU time]{\includegraphics[width=0.23\textwidth]{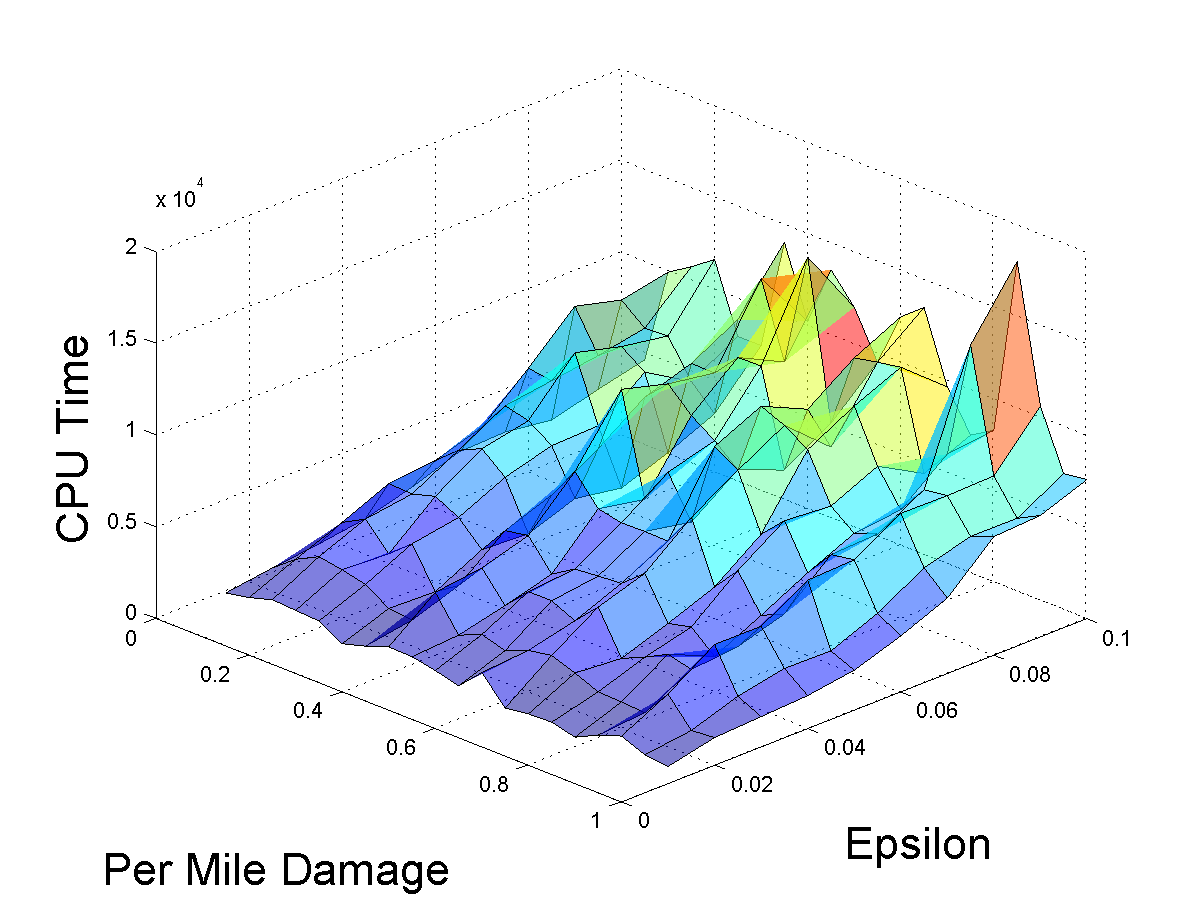}}\label{fig:12:a}
   \subfigure[Objective value]{\includegraphics[width=0.23\textwidth]{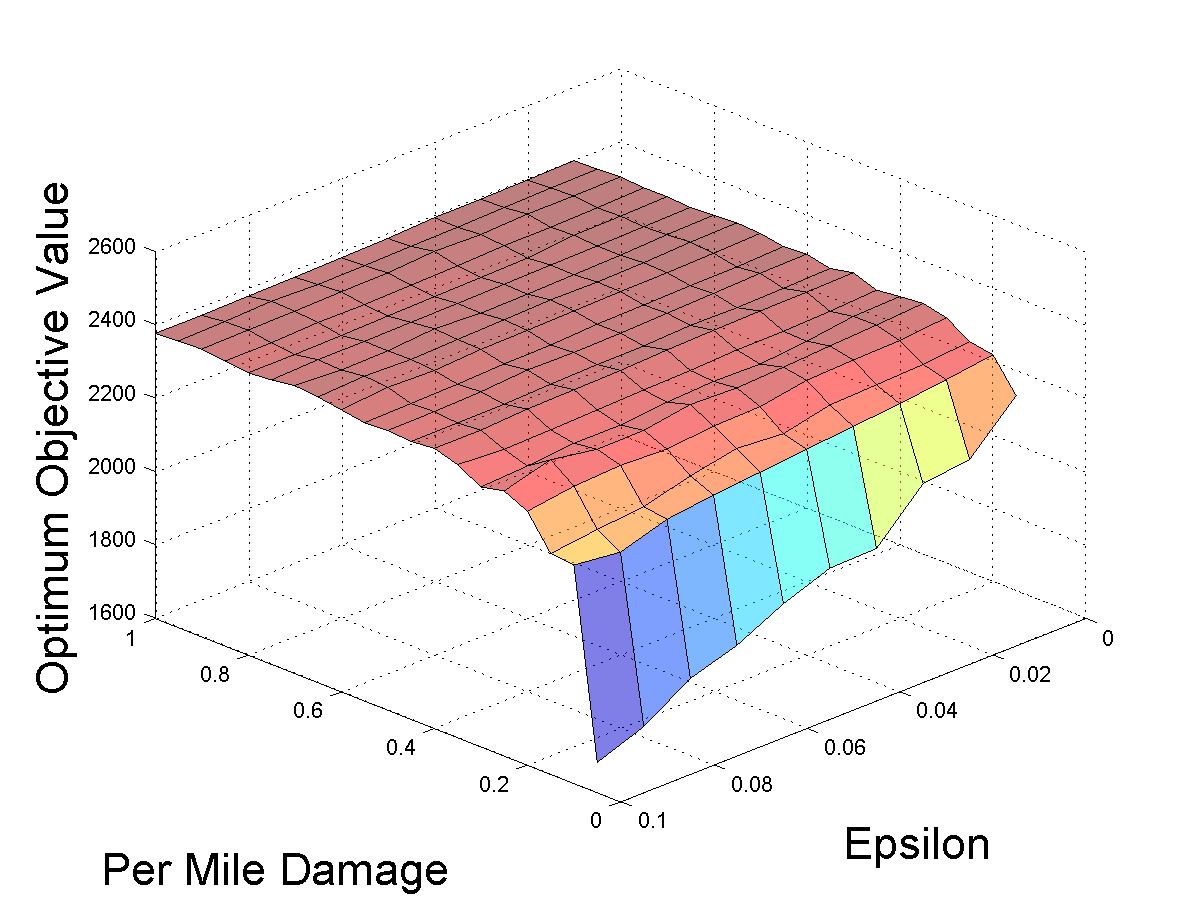}}\label{fig:12:b}
   \caption{These figures show how the CPU time and solution quality changes when chance constraints ($\epsilon$) is modified for the Rural network, when hardened lines are not damageable. These plots are generated by SBD.  }

   \label{fig:14}
\end{figure}


\section{Conclusions}

We formulated, proposed and tested new algorithms to solve the ORDGDP. Our primary contribution is an algorithm that combines the benefits of an exact method based on scenario decomposition with variable neighborhood search.  This algorithm is shown to scale well to problems that are difficult for exact methods, without sacrificing solution quality. Future directions include:
\begin{itemize}
\item Including a more accurate model of the 3-phase AC power flow equations to better exclude infeasible solutions.  Options include the DistFlow approximation in \cite{BW:89} and no-good cuts.

\item Scaling to entire city-sized distribution networks.
 \item Including a variation of the restoration problem posed by \cite{CHB:12}.
\end{itemize}

\section{ Acknowledgments}
This work was supported by the Microgrid Program of the
Office of Electricity within the U.S. Department of Energy.


\bibliographystyle{aaai}
\bibliography{Grid_Resilience}

\begin{thebibliography}{}

\bibitem[\protect\citeauthoryear{Baran and Wu}{1989}]{BW:89}
Baran, M., and Wu, F.
\newblock 1989.
\newblock {Optimal Capacitor Placement on Radial Distribution Systems}.
\newblock {\em IEEE Transactions on Power Delivery} 4(1):725--734.

\bibitem[\protect\citeauthoryear{Bent, Berscheid, and Toole}{2010}]{BBT:10}
Bent, R.; Berscheid, A.; and Toole, G.~L.
\newblock 2010.
\newblock Transmission network expansion planning with simulation optimization.
\newblock In {\em AAAI}.
\newblock AAAI.

\bibitem[\protect\citeauthoryear{Chen and Phillips}{2013}]{Chen:13}
Chen, R. L.-Y., and Phillips, C.~A.
\newblock 2013.
\newblock k-edge failure resilient network design.
\newblock {\em Electronic Notes in Discrete Mathematics} 41(0):375 -- 382.

\bibitem[\protect\citeauthoryear{Chen \bgroup et al\mbox.\egroup
  }{2014}]{Chen2014a}
Chen, R. L.-Y.; Cohn, A.; Fan, N.; and Pinar, A.
\newblock 2014.
\newblock {Contingency-Risk Informed Power System Design}.
\newblock {\em IEEE Transactions on Power Systems} 29(5):2087--2096.

\bibitem[\protect\citeauthoryear{Coffrin, Hentenryck, and Bent}{2012}]{CHB:12}
Coffrin, C.; Hentenryck, P.~V.; and Bent, R.
\newblock 2012.
\newblock Last-mile restoration for multiple interdependent infrastructures.
\newblock In {\em AAAI}.

\bibitem[\protect\citeauthoryear{Delgadillo, Arroyo, and
  Alguacil}{2010}]{Delgadillo2010}
Delgadillo, A.; Arroyo, J.; and Alguacil, N.
\newblock 2010.
\newblock {Analysis of Electric Grid Interdiction with Line Switching}.
\newblock {\em IEEE Transactions on Power Systems} 25(2):633--641.

\bibitem[\protect\citeauthoryear{{Executive Office of the
  President}}{2013}]{WhiteHouse2013}
{Executive Office of the President}.
\newblock 2013.
\newblock Economic benefits of increasing electric grid resilience to weather
  outages.
\newblock Technical report, Executive Office of the President.

\bibitem[\protect\citeauthoryear{Garcia \bgroup et al\mbox.\egroup
  }{2000}]{Garcia2000}
Garcia, P.; Pereira, J.; Carneiro, S.; da~Costa, V.; and Martins, N.
\newblock 2000.
\newblock {Three-phase power flow calculations using the current injection
  method}.
\newblock {\em IEEE Transactions on Power Systems} 15(2):508--514.

\bibitem[\protect\citeauthoryear{Garg and Smith}{2008}]{GS:08}
Garg, M., and Smith, J.~C.
\newblock 2008.
\newblock {Models and algorithms for the design of survivable multicommodity
  flow networks with general failure scenarios}.
\newblock {\em Omega} 36(6):1057--1071.

\bibitem[\protect\citeauthoryear{Garg, Jayram, and
  Narayanaswamy}{2013}]{GJN:13}
Garg, V.~K.; Jayram, T.~S.; and Narayanaswamy, B.
\newblock 2013.
\newblock Online optimization with dynamic temporal uncertainty: Incorporating
  short term predictions for renewable integration in intelligent energy
  systems.
\newblock In {\em AAAI}.

\bibitem[\protect\citeauthoryear{Golari, Fan, and Wang}{2014}]{GFW:14}
Golari, M.; Fan, N.; and Wang, J.
\newblock 2014.
\newblock Two-stage stochastic optimal islanding operations under severe
  multiple contingencies in power grids.
\newblock {\em Electric Power Systems Research} 114(0):68 -- 77.

\bibitem[\protect\citeauthoryear{Hentenryck, Gillani, and
  Coffrin}{2012}]{HGC:12}
Hentenryck, P.~V.; Gillani, N.; and Coffrin, C.
\newblock 2012.
\newblock Joint assessment and restoration of power systems.
\newblock In {\em ECAI},  792--797.

\bibitem[\protect\citeauthoryear{Jabr}{2013}]{Jabr2013}
Jabr, R.~A.
\newblock 2013.
\newblock {Robust Transmission Network Expansion Planning With Uncertain
  Renewable Generation and Loads}.
\newblock {\em IEEE Transactions on Power Systems} 28(4):4558--4567.

\bibitem[\protect\citeauthoryear{Jain, Narayanaswamy, and
  Narahari}{2014}]{JNN:14}
Jain, S.; Narayanaswamy, B.; and Narahari, Y.
\newblock 2014.
\newblock A multiarmed bandit incentive mechanism for crowdsourcing demand
  response in smart grids.
\newblock In {\em AAAI},  721--727.

\bibitem[\protect\citeauthoryear{Johnson, Lenstra, and Kan}{1978}]{Johnson1978}
Johnson, D.~S.; Lenstra, J.~K.; and Kan, A. H. G.~R.
\newblock 1978.
\newblock {The complexity of the network design problem}.
\newblock {\em Networks} 8(4):279--285.

\bibitem[\protect\citeauthoryear{Kersting}{1991}]{Kersting1991}
Kersting, W.
\newblock 1991.
\newblock {Radial distribution test feeders}.
\newblock {\em IEEE Transactions on Power Systems} 6(3):975--985.

\bibitem[\protect\citeauthoryear{Khushalani, Solanki, and
  Schulz}{2007}]{Khushalani2007}
Khushalani, S.; Solanki, J.~M.; and Schulz, N.~N.
\newblock 2007.
\newblock Optimized restoration of unbalanced distribution systems.
\newblock {\em IEEE Transactions on Power Systems} 22(2):624--630.

\bibitem[\protect\citeauthoryear{Lazic \bgroup et al\mbox.\egroup
  }{2010}]{LHMU:10}
Lazic, J.; Hanafi, S.; Mladenovic, N.; and Urosevic, D.
\newblock 2010.
\newblock Variable neighbourhood decomposition search for 0-1 mixed integer
  programs.
\newblock {\em Comp. {\&} OR} 37(6):1055--1067.

\bibitem[\protect\citeauthoryear{Li \bgroup et al\mbox.\egroup }{2014}]{Li2014}
Li, J.; Ma, X.-Y.; Liu, C.-C.; and Schneider, K.~P.
\newblock 2014.
\newblock {Distribution System Restoration With Microgrids Using Spanning Tree
  Search}.
\newblock {\em IEEE Transactions on Power Systems} PP(99):1--9.

\bibitem[\protect\citeauthoryear{Mansfield and Linzey}{2013}]{Sandy:13}
Mansfield, M., and Linzey, W.
\newblock 2013.
\newblock Hurricane sandy multi-state outage \& restoration report.
\newblock Technical Report 9308, National Association of State Eenergy
  Officials.

\bibitem[\protect\citeauthoryear{Munoz \bgroup et al\mbox.\egroup
  }{2014}]{Munoz2014}
Munoz, F.~D.; Hobbs, B.~F.; Ho, J.~L.; and Kasina, S.
\newblock 2014.
\newblock {An Engineering-Economic Approach to Transmission Planning Under
  Market and Regulatory Uncertainties: WECC Case Study}.
\newblock {\em IEEE Transactions on Power Systems} 29(1):307--317.

\bibitem[\protect\citeauthoryear{Nace \bgroup et al\mbox.\egroup
  }{2013}]{NPTZ:13}
Nace, D.; Pi�ro, M.; Tomaszewski, A.; and Zotkiewicz, M.
\newblock 2013.
\newblock Complexity of a classical flow restoration problem.
\newblock {\em Networks} 62(2):149--160.

\bibitem[\protect\citeauthoryear{Raidl, Baumhauer, and Hu}{2014}]{RBH:14}
Raidl, G.~R.; Baumhauer, T.; and Hu, B.
\newblock 2014.
\newblock Speeding up logic-based benders' decomposition by a metaheuristic for
  a bi-level capacitated vehicle routing problem.
\newblock In {\em {HM} 2014},  183--197.

\bibitem[\protect\citeauthoryear{Reddy and Veloso}{2011}]{RV:11}
Reddy, P.~P., and Veloso, M.~M.
\newblock 2011.
\newblock Strategy learning for autonomous agents in smart grid markets.
\newblock In {\em IJCAI},  1446--1451.

\bibitem[\protect\citeauthoryear{Reddy and Veloso}{2012}]{RV:12}
Reddy, P.~P., and Veloso, M.~M.
\newblock 2012.
\newblock Factored models for multiscale decision-making in smart grid
  customers.
\newblock In {\em AAAI}.

\bibitem[\protect\citeauthoryear{Reddy and Veloso}{2013}]{RV:13}
Reddy, P.~P., and Veloso, M.~M.
\newblock 2013.
\newblock Negotiated learning for smart grid agents: Entity selection based on
  dynamic partially observable features.
\newblock In {\em AAAI}.

\bibitem[\protect\citeauthoryear{Sa}{2002}]{Sa_Thesis}
Sa, Y.
\newblock 2002.
\newblock {\em Reliability Analysis of Electric Distribution Lines}.
\newblock Ph.D. Dissertation, McGill University, Montreal, Canada.

\bibitem[\protect\citeauthoryear{Salmeron, Wood, and
  Baldick}{2009}]{Salmeron2009}
Salmeron, J.; Wood, K.; and Baldick, R.
\newblock 2009.
\newblock {Worst-case interdiction analysis of large-scale electric power
  grids}.
\newblock {\em IEEE Transactions on Power Systems} 24(1):96--104.

\bibitem[\protect\citeauthoryear{Santoso \bgroup et al\mbox.\egroup
  }{2003}]{SAGS:03}
Santoso, T.; Ahmed, S.; Goetschalckx, M.; and Shapiro, A.
\newblock 2003.
\newblock A stochastic programming approach for supply chain network design
  under uncertainty.
\newblock Stochastic Programming E-Print Series. Institut für Mathematik.

\bibitem[\protect\citeauthoryear{Shann and Seuken}{2013}]{SS:13}
Shann, M., and Seuken, S.
\newblock 2013.
\newblock An active learning approach to home heating in the smart grid.
\newblock In {\em IJCAI}.

\bibitem[\protect\citeauthoryear{Thi�baux \bgroup et al\mbox.\egroup
  }{2013}]{TCHS:13}
Thi�baux, S.; Coffrin, C.; Hijazi, H.; and Slaney, J.~K.
\newblock 2013.
\newblock Planning with mip for supply restoration in power distribution
  systems.
\newblock In {\em IJCAI}.

\bibitem[\protect\citeauthoryear{Tomaszewski, Pi\'{o}ro, and
  \.{Z}otkiewicz}{2010}]{TPZ:10}
Tomaszewski, A.; Pi\'{o}ro, M.; and \.{Z}otkiewicz, M.
\newblock 2010.
\newblock {On the complexity of resilient network design}.
\newblock {\em Networks} 55(2):108--118.

\bibitem[\protect\citeauthoryear{{US Department of
  Energy}}{2013}]{DOEResilience2013}
{US Department of Energy}.
\newblock 2013.
\newblock {U.S. E}nergy sector vulnerabilities to climate change and extreme
  weather.
\newblock Technical Report DOE/PI-0013, US Department of Energy.

\bibitem[\protect\citeauthoryear{Vanderbeck and Wolsey}{2010}]{VW:10}
Vanderbeck, F., and Wolsey, L.~A.
\newblock 2010.
\newblock Reformulation and decomposition of integer programs.
\newblock In {\em 50 Years of Integer Programming}.
\newblock  431--502.

\end{thebibliography}

\end{document}